\documentclass[12pt, leqno]{amsart}  
\usepackage{amsmath,amstext,amsthm,amssymb,amsxtra}
\usepackage[top=1.5in, bottom=1.5in, left=1in, right=1in]	{geometry}
\usepackage{txfonts} % pxfonts txfonts 
\usepackage[T1]{fontenc}
\usepackage{lmodern}

 \usepackage{euler}   % better than the option below
\usepackage{pdfsync}

\usepackage{float} %To specify placement of figures

\usepackage{tikz}

\usepackage{mathtools}
\mathtoolsset{showonlyrefs,showmanualtags}

\usepackage{hyperref} %,hypertexnames=false,colorlinks,[pagebackref]
\hypersetup{
    colorlinks=true,       % false: boxed links; true: colored links
    linkcolor=blue,          % color of internal links
    citecolor=magenta,        % color of links to bibliography
    filecolor=magenta,      % color of file links
    urlcolor=cyan           % color of external links
}

\usepackage[msc-links]{amsrefs}  
\newcommand{\norm}[1]{\ensuremath{\left\|#1\right\|}}
\newcommand{\abs}[1]{\ensuremath{\left\vert#1\right\vert}}

\newcommand{\avg}[1]{\langle #1 \rangle}

\newcommand{\eps}{\varepsilon}
\newcommand{\ad}{\frac{\alpha}{d}}

\newcommand{\unit}{1\!\!1}

\newcommand{\ars}{\alpha, \mathcal{S}}

{\end{list}}

\numberwithin{equation}{section}

\newtheorem{thm}{Theorem}[section]
\newtheorem{lm}[thm]{Lemma}

\newtheorem{conj}[thm]{Conjecture}

\newtheorem{prop}[thm]{Proposition}
\newtheorem*{prop*}{Proposition}

\theoremstyle{remark}

\newtheorem*{rem*}{Remark}

 \theoremstyle{remark}
 
\numberwithin{equation}{section}

\title[Off--Diagonal Bumps] {Off--Diagonal Two Weight Bumps for Fractional Sparse Operators}

\author{Rob Rahm}
\address{Texas A\&M Mathematics Department}
\email{robrahm@math.tamu.edu}

\begin{document}

\begin{abstract}
In this paper, we continue some recent work on two weight boundedness of 
sparse operators to the "off--diagonal" setting. We use the new "entropy bumps" 
introduced in by Treil--Volberg (\cite{TreVol2016}) and improved by Lacey--Spencer 
(\cite{LacSpe2015}) and the "direct comparison bumps" introduced by Rahm--Spencer 
(\cite{RahSpe2018}) and improved by Lerner (\cite{Ler2020}). Our results are 
"sharp" in the sense that they are sharp in various particular cases. 
A feature is that given the current machinery and advances, the proofs are almost 
trivial.
\end{abstract}
\maketitle
\section{Introduction}
The topic of this article is two--weight bump conditions for sparse operators
in the ``off--diagonal'' setting (i.e $q > p$). 
We continue the line of investigation 
concerning entropy bumps that began with Treil--Volberg in \cite{TreVol2016} 
and was continued in \cites{LacSpe2015, RahSpe2016, RahSpe2018} and the 
line of investigation concerning "direct comparison bumps" introduced in 
\cite{RahSpe2018} and continued and improved in \cite{Ler2020}. 

One innovation that needs to be mentioned is that all of the results 
apply to general measures, rather than to measures that are absolutely 
continuous and have a density. Throughout, $\mu$ and $\nu$ are two locally finite 
measures. Quantities like 
$\avg{\nu f}_{Q}$ mean $\frac{1}{\abs{Q}}\int_{Q}f d\nu$. And the maximal function 
applied to $\nu$ (i.e. $M\nu(x)$ is defined using the averages 
$\avg{\nu}_{Q} := \frac{1}{\abs{Q}}\int_{Q}d\nu$). 

We will be concerned with sparse 
operators of the form ( $0\leq \alpha < d$): 
\begin{align*}
T_{\ars}\nu f
:= \sum_{Q\in\mathcal{S}}(\abs{Q}^{\frac{\alpha}{d}}\avg{\abs{\nu f}}_{Q})
    \unit_{Q}, 
\hspace{.5in}
\textnormal{ where }
\hspace{.5in}
\avg{\nu f}_{Q}:= \frac{1}{\abs{Q}}\int_{Q}fd\nu.
\end{align*}

A more general class of these operators are studied in \cite{FacHyt2018}. 
Unlike that paper, we make no 
assumption about the weights being in $A_\infty$ and we only concentrate 
on the off--diagonal case.

There are many results like this in the literature. 
The main contributions of this paper 
are (1) the results are new; (2) the proofs are succinct and highlight the development 
of this area and (3) indicate the challenges in proving various "separated bump 
theorems".

The first two deal with entropy bumps. It was proven in the $p=q=2$ case in 
\cite{TreVol2016} and in the general 
$p=q$ case in \cite{LacSpe2015}; both with $\alpha = 0$. There is also a 
version in \cite{RahSpe2016} that has the "wrong" homogeneity. 

\begin{thm}\label{T:EB0}
Let $T_{\ars}$ be a sparse operator and 
$1< p < q < \infty$ then 
\begin{align*}
\norm{T_{\alpha,\mathcal{S}}\nu\cdot:L^p(\mu)\to L^q(\nu)} 
\lesssim [\nu,\mu]_{A_{p,q},\mathcal{S}}[\nu]_{A_\infty,\mathcal{S}}^{\frac{1}{q}}
+ [\mu,\nu]_{A_{q',p'},\mathcal{S}}[\mu]_{A_\infty, \mathcal{S}}^{\frac{1}{p'}},
\end{align*}
where
\begin{align*}
[\nu,\mu]_{A_{p,q},\mathcal{S}}
:= \sup_{Q\in\mathcal{S}}
    \frac{\mu(Q)^{\frac{1}{q}}\nu(Q)^{\frac{1}{p'}}}{\abs{Q}^{1-\ad}}
\end{align*}  
and 
\begin{align*}
\rho(Q; \nu) := \frac{1}{\nu(Q)}\int_{P}M(\nu\unit_{P})
\hspace{.5in}
\textnormal{and}
\hspace{.5in}
[\nu]_{A_\infty,\mathcal{S}} := \sup_{Q\in\mathcal{S}}\rho(Q;\nu),
\end{align*}
\end{thm}

\begin{thm}\label{T:EB}
Let $T_{\ars}$ be a sparse operator and 
$1< p < q < \infty$ then: 
\begin{align*}
\norm{T_{\alpha,\mathcal{S}}\nu\cdot:L^p(\mu)\to L^p(\nu)} 
\lesssim S_{\mathcal{E}}^{\frac{1}{q}}\mathcal{E}_{p,q}(\mu,\nu)
    + S_{\mathcal{E}}^{\frac{1}{p'}}\mathcal{E}_{q'p'}(\nu,\mu),
\end{align*}
where: 
\begin{align*}
\mathcal{E}_{p,q}(\mu,\nu)
:= \sup_{Q}\frac{\mu(Q)^{\frac{1}{q}}\nu(Q)^{\frac{1}{p'}}}{\abs{Q}^{1-\ad}}
    \rho(Q;\nu)^{\frac{1}{q}}\eps(\rho(Q;\nu))^{\frac{1}{q}}
\end{align*}  
%and $\rho(Q; \nu)$ is the "local $A_\infty$ characteristic": 
%\begin{align*}
%\rho(Q; \nu) := \frac{1}{\nu(Q)}\int_{P}M(\nu\unit_{P})
%\hspace{.5in}
%\textnormal{and}
%\hspace{.5in}
%\rho(Q; w) := \frac{1}{\mu(Q)}\int_{P}M(\mu\unit_{P}), 
%\end{align*}
and $\eps$ is a monotonic increasing function that satisfies 
$S_{\mathcal{E}}:=\sum_{r=0}^{\infty}\eps(2^r)^{-1} <\infty.$
\end{thm}

The next theorem was introduced in \cite{RahSpe2018} and improved in 
\cite{Ler2020} (for $p=q$ and $\alpha =0$): 
\begin{thm}\label{T:CB}
Let $T_{\ars}$ be a sparse operator and 
$1 < p < q < \infty$ and $p>1$ then: 
\begin{align*}
\norm{T\nu\cdot:L^p(\mu)\to L^p(\nu)}
\lesssim 
S_{\mathcal{D}}^{\frac{1}{q}}\mathcal{D}_{p,q}(\mu,\nu) + 
S_{\mathcal{D}}^{\frac{1}{p'}}\mathcal{D}_{q',p'}(\nu,\mu),
\end{align*}
where: 
\begin{align*}
\mathcal{D}_{p,q}(\mu,\nu)
:= \sup\frac{{\mu(Q)^{\frac{1}{q}}}
    \nu(Q)^{\frac{1}{p'}}}
    {\abs{Q}^{1-\ad}}
    \eps(\avg{\nu}_{Q})^{\frac{1}{q}}.
\end{align*}

and $\eps$ is a function that is decreasing on $(0,1)$ increasing on 
$(1,\infty)$ and satisfies 
$S_{\mathcal{D}}:=\sum_{r=-\infty}^{\infty}\eps(2^r)^{-1}$. 
\end{thm}

\section{Background and Discussion }
It was noted in \cites{Muc1972,HunMucWhe1973,MucWhe1974} that the standard Muckenhoupt condition 
was necessary but not sufficient for the two--weight boundedness of operators of interest. The 
purpose of the various "bumps" are to replace the standard 
$A_p$ (or $A_{p,q}$ in our case) condition with a slightly bigger condition that is 
sufficient to give two--weight boundedness. This is a well--developed area, 
see \cites{Neu1983,Per1994,CruRezVol2014, Lac2016, NazRezTreVol2013, Ler2013, 
CruMoe2013a, CruMoe2013b, Cru2017, CruMarPer2016} and the references therein for 
more information.

The "entropy bumps" in Theorem \ref{T:EB} were introduced by Treil--Volberg in 
\cite{TreVol2016} in the $p=q=2$ case. This was extended to the general $p=q$ case 
by Lacey--Spencer in \cite{LacSpe2015}. 
These operators were also studied in the off--diagonal $p\leq q$, $\alpha>0$ setting 
in \cites{RahSpe2016,PanSun2019}.
Our theorem here is "sharp" in the sense that when 
$\nu$ and $\mu$ are $A_\infty$, we recover the 
sharp results of, for example \cites{CruMoe2013a, CruMoe2013b, LacMoePerTor2010}. 
The "direct comparison bumps" were introduced in \cite{RahSpe2018} and improved in 
\cite{Ler2020}. 

One feature of the proofs is that they \textit{only} require $A_{p,q}$ and 
$A_\infty$ data for the cubes in the sparse collection. Compare this with the 
proofs in, for example, \cite{FacHyt2018} that require the weight be in 
$A_\infty$. This small difference is what allows us to do the "bootstrapping" 
argument in the proof of Theorem \ref{T:EB}.

In addition, it is
the off--diagonal setting that allows us to replace the normal norm in the testing 
inequalities with $L^1$ norms (this is Proposition \ref{P:T} below; see also 
\cite{CruMoe2013b}*{Theorem 1.1} and similar ideas in \cite{CruMarPer2016}). This is why off--diagonal 
results are sharper than 
on--diagonal results. Indeed, in the proof of Theorem \ref{T:EB}, 
if we had to work with $L^q$ norms in the testing constants, we would have to use the 
triangle inequality 
to estimate the sum over $a\geq 0$ and this would require that $\sum_{a\geq 0}
\eps(2^a)^{-\frac{1}{q}}$ be finite -- and this is a stronger assumption than 
what we have here (indeed: this is exactly what happens in \cite{RahSpe2016} 
in the on--diagonal setting).
Based on comparisons with Orlicz conjectures, the following conjecture is made 
(Lerner almost proves this  in \cite{Ler2020}):

\begin{conj}\label{C:DC}
$\norm{T\nu\cdot:L^p(\nu)\to L^p(\mu)} 
\lesssim \mathcal{D}_{p,p}(\nu,\mu) 
    + \mathcal{D}^{\ast}_{p',p'}(\mu,\nu)$.
\end{conj}

The Orlicz, entropy, and direct comparison bumps are not strictly comparable. The 
Orlicz bumps are  the most established but require the most in terms of local 
integrability (and they apply only to measures with densities).  The entropy bumps 
are guaranteed to be bounded in the one weight setting (and record important $A_\infty$
information about the operator norms).  The direct comparison bumps require the least 
in terms of local integrabilty but are not as popular and the information 
they record (comparison to Lebesgue measure) does not seem to be as relevant 
as $A_\infty$ data in the weighted theory. An advantage is that they are easier 
to verify than Orlicz bumps or entropy bumps. The measures
\begin{align*}
d\nu(x) = \frac{dx}{x(1-\log x)^2},
\hspace{.5in}
d\mu(x) = x^2dx
\end{align*}
satisfy the two weight conditions in \cite{Ler2020} or even in \cite{RahSpe2018}. Yet 
$\nu$ is not in $L\log L$ and so neither the Orlicz nor the entropy bumps will detect 
the boundedness of sparse operators with these weights.

\textbf{Acknowledgment.} I'd like to thank David Cruz-Uribe for some comments about this
paper (in particular, the observation that it applies to measures and not just absolutely 
continuous measures).

\section{Preliminaries and Notation}\label{S:PL}
A collection of cubes, $\mathcal{S}$, is called $\lambda$--Sparse ($0<\lambda<1$) if 
for every $Q\in\mathcal{S}$ there is a set $E_Q\subset Q$ with 
$\abs{E_Q} \geq \lambda\abs{Q}$ and the sets $\{E_Q:Q\in\mathcal{S}\}$ are pairwise 
disjoint. If $R\in\mathcal{S}$ then using 
 the pairwise disjointness of the $\{E_Q\}$ we have the following well--know estimate:
\begin{align}\label{E:r}
\sum_{Q\in\mathcal{Q}:Q\subset R}\nu(Q)
%\simeq \sum_{Q\in\mathcal{Q}:Q\subset Q_0}\avg{\nu}_{Q} \abs{E_Q}
\simeq \int_{R}\sum_{Q\in\mathcal{Q}:Q\subset R}\avg{\nu}_{Q}\unit_{E_Q}
\leq \int_{R}M(\nu \unit_{R})
= \rho(R;\nu)\nu(R).
\end{align}

A consequence of \cite{CruMoe2013b}*{Theorem 1.1} is (note the $L^1$ norms): 
\begin{prop}\label{P:T}
For $1 < p < q$ there holds
\begin{align*}
\norm{T_{\ars}(\nu \cdot):L^p(\nu)\to L^q(\mu)}
\lesssim 
\mathcal{T} + \mathcal{T}^\ast
\end{align*}
where: 
\begin{align*}
\mathcal{T} 
&:= \sup_{\mathcal{S} \textnormal{ is sparse}}\sup_{R\in\mathcal{S}}\nu(R)^{-\frac{1}{p}}
  \norm{\sum_{Q\in\mathcal{S}:Q\subset R}(\abs{Q}^{\ad}\avg{\nu\unit_{R}}_{Q})^{q}\unit_{E_Q}}_{L^{1}(\mu)}^{\frac{1}{q}}
\\\mathcal{T}^{\ast}
&= \sup_{\mathcal{S} \textnormal{ is sparse}}\sup_{R\in\mathcal{S}}w(R)^{-1/q'}
  \norm{\sum_{Q\in\mathcal{S}:Q\subset R}
  (\abs{Q}^{\ad}\avg{\mu}_{Q})^{p'}\unit_{E_Q}}_{L^{1}(\nu)}^{\frac{1}{p'}}. 
\end{align*}
\end{prop}

\begin{proof}
Indeed, \cite{CruMoe2013b}*{Theorem 1.1} says that if $p < q$ then:
\begin{align*}
\norm{T_{{\ars}}\nu\cdot: L^p(\nu) \to L^q(\mu)}
\simeq 
\norm{M_{{\alpha}}\nu\cdot: L^p(\nu) \to L^q(\mu)} + 
\norm{M_{{\alpha}}\nu\cdot: L^{q'}(\mu) \to L^{p'}(\nu)},
\end{align*}
where 
$M_\alpha f (x) := \sup_{Q}\abs{Q}^{\frac{\alpha}{d}}\avg{\abs{f}}_{Q}\unit_{Q}(x)$
is the fractional maximal operator. The well--known (\cite{Saw1982}) testing 
conditions for this operator reduce to the ones in Proposition \ref{P:T}. This 
is because the linearization of $M_\alpha$ is the function inside the norm 
in this proposition and the pairwise 
disjointness of the $\{E_Q\}$ allows us to pass the exponent $q$ (and $p'$) 
under the sum.
\end{proof}
\section{Proofs of Theorems \ref{T:EB0} and \ref{T:EB}}\label{S:EB0}
We begin with a lemma that will be used in the proof of both of 
these theorems.
\begin{lm}\label{L:ml}
Let $\mathcal{S}$ be a sparse collection. With notation as established 
above, for every $R\in\mathcal{S}$ there holds: 
\begin{align*}
\sum_{Q\in\mathcal{S}:Q\subset R}\abs{Q}^{q\ad}\avg{\nu}_{Q}^{q}\mu(Q)
\lesssim [\nu,\mu]_{p,q,\mathcal{S}}^{q}
    [\nu]_{A_\infty, \mathcal{S}}\nu(R)^{\frac{q}{p}}.
\end{align*}
\end{lm}
\begin{proof}[Proof of Lemma \ref{L:ml}]
Let $\mathcal{S}^\ast$ be the maximal cubes in $\mathcal{S}$ contained 
in $R$. 
The sum above can be organized as follows: 
\begin{align*}
\sum_{Q\in\mathcal{S}:Q\subset R}\abs{Q}^{q\ad}\avg{\nu}_{Q}^{q}\mu(Q)
=\sum_{Q^\ast\in\mathcal{S}^\ast}
\sum_{Q\in\mathcal{S}:Q\subset Q^\ast}(\frac{\mu(Q)\nu(Q)^{\frac{q}{p'}}}
            {\abs{Q}^{q - q\ad}})\nu(Q)^{\frac{q}{p}}.
\end{align*}
The inner sum is dominated by:
\begin{align*}
[\nu,\mu]_{p,q,\mathcal{S}}^{q}
\sum_{Q\in\mathcal{S}_a:Q\subset Q^\ast}\nu(Q)^{\frac{q}{p}}
\leq [\nu,\mu]_{p,q,\mathcal{S}}^{q} \nu(Q^\ast)^{\frac{q}{p}-1}
     \sum_{Q\in\mathcal{S}_a:Q\subset Q^\ast}\nu(Q)
 \lesssim [\nu,\mu]_{p,q,\mathcal{S}}^{q}[\nu]_{A_\infty, 
 \mathcal{S}}\nu(Q^\ast)^{\frac{q}{p}}
\end{align*}
The "$\leq$" is trivial and the "$\lesssim$" uses \eqref{E:r}. 
Using the fact that $\frac{q}{p}> 1$. Using the maximality of the cubes $Q^\ast$ 
this can be summed over $Q^\ast$ in $\mathcal{S}^\ast$ to 
$[\nu,\mu]_{p,q,\mathcal{S}}^{q}[\nu]_{A_\infty, \mathcal{S}}\nu(R)^{\frac{q}{p}}$. 
\end{proof}

\begin{proof}[Proof of Theorem \ref{T:EB0}]
We show 
$\mathcal{T}\lesssim[\nu,\mu]_{p,q,\mathcal{S}}
[\nu]_{A_\infty,\mathcal{S}}^{\frac{1}{q}}$ 
(the estimate $\mathcal{T}^\ast \lesssim 
[\mu,\nu]_{q',p',\mathcal{S}}[\mu]_{A_\infty,\mathcal{S}}^{\frac{1}{p'}}$ is dual). 
Observe that we can estimate the $q^{\textnormal{th}}$ power of the norm in the 
definition of $\mathcal{T}$ using Lemma \ref{L:ml} as
\begin{align*}
\sum_{Q\in\mathcal{S}:Q\subset R}\abs{Q}^{q\ad}\avg{\nu}_{Q}^{q}\mu(Q)
\lesssim[\nu,\mu]_{p,q,\mathcal{S}}^{q}[\nu]_{A_\infty, 
\mathcal{S}}\nu(R)^{\frac{q}{p}}.
\end{align*}
Taking $q^{\textnormal{th}}$ roots we conclude that 
$\mathcal{T}\lesssim [\nu,\mu]_{p,q,\mathcal{S}}
[\nu]_{A_\infty, \mathcal{S}}^{\frac{1}{q}}$ as 
claimed. 
\end{proof}

\begin{proof}[Proof of Theorem \ref{T:EB}]
$\mathcal{T}\lesssim\mathcal{E}_{p,q}(\nu,\mu)$ 
(the estimate $\mathcal{T}^\ast \lesssim 
\mathcal{E}_{q',p'}(\mu,\nu)$ is dual). 
Let $\mathcal{S}_a$ be those cubes with $2^a < \rho(Q;\nu)\leq 2^{a+1}$.
Using Lemma \ref{L:ml}, the 
$q^{\textnormal{th}}$ power of the norm in the 
definition of $\mathcal{T}$ can be estimated as follows
\begin{align}\label{E:EB}
\sum_{Q\in\mathcal{S}:Q\subset R}\abs{Q}^{q\ad}\avg{\nu}_{Q}^{q}\mu(Q)
=\sum_{a\geq 0}
\sum_{Q\in\mathcal{S}_a:Q\subset R}\abs{Q}^{q\ad}\avg{\nu}_{Q}^{q}\mu(Q)
\lesssim \sum_{a\geq 0}[\nu,\mu]_{p,q,\mathcal{S}_a}^{q}
    [\nu]_{A_\infty, \mathcal{S}_a}\nu(R)^{\frac{q}{p}}
\end{align}
Multiplying and dividing the summands by $\eps(2^a)$ and using 
the fact that for the cubes in question, $\eps(2^a)\leq\eps(\rho(Q;\nu))$ and 
$\rho(Q;\nu)\simeq 2^a$, the summands are estimated as:
\begin{align*}
\frac{1}{\eps(2^a)}  
    \sup_{Q\in\mathcal{S}_a}
    \left(\frac{\mu(Q)^{\frac{1}{q}}\nu(Q)^{\frac{1}{p'}}}{\abs{Q}^{1-\ad}}
    \rho(Q;\nu)^{\frac{1}{q}}\eps(\rho(Q;\nu))^{\frac{1}{q}}\right)^q 
    \nu(R)^{\frac{q}{p}}
\leq \frac{1}{\eps(2^a)}\mathcal{E}_{p,q}(\nu,\mu)^{q}\nu(R)^{\frac{q}{p}}.
\end{align*}
Using the summability of $\eps$ this can be 
summed in $a$ to $\mathcal{S}_{\mathcal{E}}\mathcal{E}_{p,q}(\nu,\mu)^{q}
\nu(R)^{\frac{q}{p}}$. Taking $q^{\textnormal{th}}$ roots gives the desired estimate.
\end{proof}

\section{Proof of Theorem \ref{T:CB}}
We will show 
$\mathcal{T}\lesssim\mathcal{D}_{p,q}(\nu,\mu)$ (the estimate 
$\mathcal{T}^\ast \lesssim \mathcal{D}_{q',p'}(\mu.\nu)$ is dual). 
Let $\mathcal{S}_a$ be those cubes with $2^{a}<\avg{\nu}_{Q}\leq 2^{a+1}$. 
(observe that $-\infty < r< \infty$) and let $\mathcal{S}_a^\ast$ be the maximal cubes 
in $\mathcal{S}_a$. As above, the $q^{\textnormal{th}}$ power of the
norm in the definition of $\mathcal{T}$ can be 
organized as follows: 
\begin{align}\label{E:DB}
\sum_{Q\in\mathcal{S}:Q\subset R}\abs{Q}^{q\ad}\avg{\nu}_{Q}^{q}\mu(Q)
=\sum_{a\in \mathbb{Z}}\sum_{Q^\ast\in\mathcal{S}_a^\ast}
\sum_{Q\in\mathcal{S}_a:Q\subset Q^\ast}\abs{Q}^{q\ad}\avg{\nu}_{Q}^{q}\mu(Q).
\end{align}
Concerning the inner sum, this is (similar to the above):
\begin{align*}
\sum_{Q\in\mathcal{S}_a:Q\subset Q^\ast}(
            \frac{\mu(Q)\nu(Q)^{\frac{q}{p'}}}
            {\abs{Q}^{q - q\ad}})\nu(Q)^{\frac{q}{p}}
%&\lesssim \mathcal{D}^{q}
%\frac{1}{\eps(2^a)}
%\sum_{Q\in\mathcal{S}_a:Q\subset Q^\ast}\nu(Q)^{\frac{q}{p}}
&\lesssim \mathcal{D}_{p,q}(\nu,\mu)^{q}\frac{1}{\eps(2^a)} \nu(Q^\ast)^{\frac{q}{p}-1}
     \sum_{Q\in\mathcal{S}_r:Q\subset Q^\ast}\nu(Q)
\\&\leq \mathcal{D}_{p,q}(\nu,\mu)^{q}\frac{\nu(Q^\ast)^{\frac{q}{p}}}{\eps(2^a)},
\end{align*}
where for the "$\lesssim$" we multiplied and divided by 
$\eps(2^a)$ and used the fact that 
$\avg{\nu}_{Q}\simeq 2^a$. The ``$\leq$'' is obtained using sparseness 
combined with $\nu(Q)\simeq 2^a\abs{Q}$ for the cubes in question.
Similar to the above proof, the summability condition on $\eps$ implies that
$\mathcal{T}\lesssim S_{\mathcal{D}}^{\frac{1}{q}}\mathcal{D}_{p,q}(\nu,\mu)$.

%%%%%%%%%%%%
%%%References%%%
%%%%%%%%%%%%

%%%%%%%%%%%%
%%%%END%%%%%%
%%%%%%%%%%%%
\end{document}